\documentclass[12pt]{article}
\usepackage{amssymb}
\usepackage{enumerate}
\usepackage{graphicx}
\usepackage{amsthm}
\usepackage{color}

\def\R{{\rm I\! R}}

\def\C{\mbox{l\hspace{-.47em}C}}

\def\dst{\displaystyle}

\newtheorem{theorem}{Theorem}
\newtheorem{corollary}{Corollary}

\newtheorem{lemma}{Lemma}

\newtheorem{remark}{Remark}

\title{Global injectivity of planar non-singular maps polynomial in one variable
\thanks{Dipartimento di Matematica, Universit\`a di Trento, I-38121 Trento (TN), Italy; email: marco.sabatini@unitn.it; MSC Classification: 14R15, 26B10.    }}

\author{Marco Sabatini  }
\begin{document}

\maketitle
\begin{abstract} We consider non-singular and Jacobian maps whose components are polynomial in the variable $y$. We  prove that if a map has $y$-degree one, then it is the composition of a triangular map and a quasi-triangular map. We also prove that non-singular $y$-quadratic maps are injective if one of the leading functional coefficients does not vanish.  Moreover, $y$-quadratic Jacobian maps are the composition of a quasi-triangular map and 3 triangular maps. Other results are given for wider classes of non-singular maps, considering also injectivity on vertical strips $I \times \R$.

\end{abstract}

\section{Introduction}

Let us consider a map $F\in C^1(\Omega,\R^2)$, $\Omega$ open connected subset of $\R^2$. Let
$$
J_F (x,y) = \left( \matrix{P_x(x,y) & P_y(x,y) \cr Q_x(x,y) & Q_y(x,y)}  \right)
$$ 
be the jacobian matrix of $F$ at $(x,y)$. We denote by $d_F(x,y) = \det J_F(x,y)  = P_x(x,y)\, Q_y(x,y) - P_y(x,y)\, Q_x(x,y)$ its determinant. 
We say that $F(x,y)$ is a \textit{non-singular map} if $d_F(x,y) \neq 0$ on all of $\Omega$.
We say that $F(x,y)$ is a \textit{Jacobian map} if $d_F(x,y) $ is a non-zero constant on all of $\Omega$. The implicit function theorem gives the  invertibility of a map in a neighbourhood of a point   $(x^*,y^*) \in \Omega$  such that  $d_F(x^*,y^*) \neq 0$. On the other hand, even if $d_F(x,y) \neq 0$ on all of $\Omega$, the map can be non-invertible, as shows the exponential map $(e^y \cos x,e^y \sin x)$. The search for additional conditions ensuring the global invertibility or injectivity of a locally invertible map is a classical problem. A fundamental result is Hadamard global inverse function theorem, which gives the global invertibility of a proper non-singular map $F\in C^1(\R^n,\R^n)$. Proving properness is often too difficult, so that different paths to the global invertibility or injectivity are looked for.
In this field some old  problems still resist the attempts to find  a solution. The celebrated Jacobian Conjecture is  concerned with polynomial maps $F: \C^n \to \C^n$ \cite{Ke}. According to such a conjecture, every polynomial map with non-zero constant jacobian determinant should be invertible, with polynomial inverse. Such a statement and its variants were studied  in several settings, even replacing $\C^n$ with $\R^n$ or other fields, and several partial results were proved, but it is not yet proved or disproved even for $n=2$ \cite{BCW,vdE}. It is known that planar polynomial Jacobian maps are invertible if the map has degree $\leq 100$, or one componente has prime degree, or one component has  degree 4. See  \cite{BCW,vdE} for other partial results.

A stronger statement known as the Real Jacobian Conjecture, in which the condition $d_F = const. $ is replaced by $d_F \neq 0$ was proved to be false \cite{Pi1}. Recent research found additional conditions under which a non-singular planar polynomial map is injective \cite{BGL,BV}. In particular it is known that non-singular planar polynomial maps of degree $\leq 4$ are injective \cite{BO}, as well as maps with one component of degree $\leq 3$ \cite{BdSF}. Such results are concerned  with the maps invertibility, not dealing with the  polynomiality of the inverse maps.

Recently some attention has been devoted to prove the existence of unbounded injectivity regions  of polynomial non-singular maps, too  \cite{MX}. Another research line is directed to more general settings, as the study of semi-algebraic maps \cite{BGV} or quasi-polynomial maps \cite{Pi2}.

In this paper we are concerned with maps polynomial in one variable but not necessarily polynomial in the other one, as in \cite{Pi2}.  
In other words, we consider maps of the form
\begin{equation}  \label{F} 
F(x,y) = \left( \sum_{i=0}^m p_i(x) y^i , \sum_{j=0}^n q_j(x) y^j \right)  ,
\end{equation}
called maps of type $(m,n)$. Differently from \cite{Pi2} we do not consider only Jacobian maps. In order to study also maps with some $p_i(x)$ or $q_j(x)$ not defined on all of $\R$, as it often occurs with rational functions, we  work on  domains obtained as cartesian products of an open interval (possibily unbounded) and the real line, $\Sigma = I \times \R$. Such domains are usually called \textit{strips}. This allows also to study injectivity regions of maps  non-singular only on subregions of their domains, as 
$F(x,y) = (x^2, x^2y)$. Such a map is non-singular for $x \neq 0$, and is invertible both for $x > 0$ and for $x < 0$, with inverse maps $\dst{ F_\pm(u,v) = \left( \pm \sqrt{u}, \frac vu  \right)  }$.

In section 1 we give some preliminary lemmata about non-singular maps, then we prove that every non-singular map of type $(1,1)$ is the composition of a triangular map and a quasi-triangular map. If additionally $F$ is polynomial and $F(0,0)=0$, then it is the composition of two triangular maps, hence it has  a polynomial inverse. After that we study the injectivity of a class of non-singular maps of type $(m,m)$. 
As a consequence, we get the injectivity of non-singular $y$-quadratic maps with non-vanishing leading coefficient. We also consider maps of type $(2hl,2h)$.

In section 2 we are concerned with maps whose Jacobian determinant does not depend on $y$, i. e.  $d_F(x,y) = \delta(x)$, called $\delta$-maps . Jacobian maps are a special case of such maps. We prove that a non-singular $\delta$-map of type $(m,1)$   is the  composition of a quasi-triangular map and $m$   triangular maps. If additionally $F(x,y)$ is polynomial and $F(0,0)=(0,0)$, then it is the composition of $m+1$   triangular maps, hence it has  a polynomial inverse. As a corollary we prove that a special class of non-singular $\delta$-maps of type $(m,m)$ is the  composition of a quasi-triangular map and $m+1$   triangular maps. If additionally  $F(x,y)$ is polynomial and $F(0,0)=(0,0)$, then it is the composition of $m+2$   triangular maps. As a corollary we have that a $y$-quadratic non-singular $\delta$-map is the composition of a quasi-triangular map and 3 triangular maps. If additionally $F(x,y)$ is polynomial and $F(0,0)=(0,0)$, then it is the composition of 4 triangular maps, hence it has a polynomial inverse. Finally we consider a class of non-singular $\delta$-map of type $(Lm,m)$.

\section{Non-singular $y$-polynomial maps}

We consider  maps $F \in C^1(\Sigma , \R^2)$, $F(x,y) = (P(x,y),Q(x,y))$, where $\Sigma = I \times \R$, $I = (I_-,I_+)$, possibly with $I_\pm = \pm \infty$. We denote partial derivatives by subscripts. 
We are concerned with maps whose components are both polynomial in one variable, that we choose to be $y$. Hence we can write
\begin{equation}  \label{F} 
F(x,y) = \left( \sum_{i=1}^m p_i(x) y^i , \sum_{j=1}^n q_j(x) y^j \right)  .
\end{equation}
with $p_i(x) \in C^1(I,\R)$ for $i=1, \dots, m$,  $q_j(x) \in C^1(I,\R)$ for $j=1, \dots, n$. Such a map is said to be of type $(m,n)$. For the sake of simplicity in the following we shall often omit the dependence on $x$ of the $p_i$ and the $q_j$. The functions   $p_i$ for $i=1, \dots, m$,  $q_j  $ for $j=1, \dots, n$ will be called \textit{coefficients} of $F(x,y)$.  
The functions $p_m$ and  $q_n$ are said to be the \textit{leading coefficients of $F(x,y)$}. 
When dealing with polynomial maps, i. e. with maps polynomial in both variables,  it will be implicitly assumed that $I = \R$, so that $\Sigma = \R^2$.

We say that  (\ref{F}) is a \textit{ quasi-triangular map} if it has the form $F(x,y) = (\alpha(x), by+\beta(x))$ or $F(x,y) = (ax +\alpha(y), \beta(y))$ for some non-zero $a,b \in \R$ and some  functions $\alpha, \beta$.  Every non-singular quasi-triangular map is injective, hence invertible on its image. Its inverse  $F^{-1}$ is itself a quasi-triangular map of the same type. For instance, if   $F(x,y) = (\alpha(x), by+\beta(x))$, its Jacobian determinant is $\alpha'(x)\, b \neq 0$, hence $\alpha(x)$ is strictly monotone,  and the inverse map is
$$
F^{-1} (u,v) = \left( \alpha^{-1}(u), \frac {v - \beta(\alpha^{-1}(u))}b   \right) ,
$$
defined on the strip $\alpha(\R) \times \R$.
Similarly if $F(x,y) = (ax +\alpha(y), \beta(y))$.
We say that  (\ref{F}) is a \textit{  triangular map} if it has the form $F(x,y) = (a x, b y+\beta(x))$ or $F(x,y) = (a x +\alpha(y), b y)$ for some non-zero $a, b \in \R$ and some functions $\alpha, \beta$. Every   triangular map is Jacobian and  invertible, having a   triangular map as inverse, defined on all of $\R^2$.

Some of next proofs use the following Lemma, based on an argument often used studying planar non-singular maps \cite{MO,S1}.

\begin{lemma}  \label{lemmauno} 
Let $F\in C^1(\Sigma,\R^2)$ be a non-singular map. If the level sets of $P(x,y)$ $\Big( Q(x,y) \Big)$ are connected, then $F(x,y)$ is injective. 
\end{lemma}
{\it Proof.}   
In order to prove that for every $(\overline {p},\overline {q}) \in \R^2$ there exists at most one point $(\overline {x},\overline {y})$ such that  $F(\overline {x},\overline {y}) = (\overline {p},\overline {q}) $, let us consider the Hamiltonian system
\begin{equation} \label {ham}
\left\{   \matrix{  \dot x = - P_y (x,y) \hfill \cr \dot y = P_x(x,y) \hfill .}  \right.
\end{equation}
It has no critical points and its orbits are contained in the level sets of $P(x,y)$, which are connected, hence every orbit coincides with a level set $P(x,y) = \overline {p}$. If $F(\overline {x},\overline {y}) = (\overline {p},\overline {q}) $, then 
$(\overline {x},\overline {y}) $ belongs to the level set $P(x,y) = \overline {p}$. 
The derivative of $Q(x,y)$ along the solutions of (\ref{ham}), i.e. $\dot Q = - Q_x P_y + Q_y P_x = d \neq 0$,   does not vanish, hence $Q(x,y)$ is strictly monotone on every  orbit, so that the equation $Q(x,y) = \overline {q}$ has at most one solution on such an orbit. This gives the injectivity of $F(x,y)$.  
\hfill  $\clubsuit$ \\

\begin{remark} \label{nonconnected}The level sets of the components of non-singular maps are not necessarily connected. The map $F(x,y) = (y^2 - e^x, (y-1)^2 - e^x)$ is non-singular and both components have disconnected zero-level sets \cite{S1}. Such a map is invertible on its image, with inverse 
$$
F^{-1} (u,v) = \left (  \frac{u-v+1}2, \ln \Big ((u-v+1)^2 - u \Big)  - \ln 4 \right).
$$
\end{remark}

\begin{lemma}  \label{lemmadue} Let $I \in \R$ be an open interval, $h,k,l $ positive integers, $p,q \in C^1(I,\R)$ such that
\begin{equation}  \label{equalemmauno} 
h p'(x) q(x) - k q'(x) p(x) = 0, \qquad x \in I.
\end{equation}
\begin{itemize}
\item[i)] If $q(x) \neq 0$ on  $I$, then there exists $c \in \R$ such that $p(x)^h = c\,  q(x)^k$   for all $x \in I$.
\item[ii)] If $p, q \in C^\omega(I,\R)$ and they do not vanish identically on $I$, then there exist $c_p, c_q \in \R$, $c_p c_q \neq 0$, such that $ c_q \, q(x) ^k + c_p\,   p(x)^h = 0$   for all $x \in I$.
\item[iii)] In both $i)$ and $ii)$, if $k = hl$, then there exist $c_l \in \R$, $c_l \neq 0$, such that $ p(x) = c_l q(x)^l$   for all $x \in I$.
\end{itemize}
\end{lemma} 
{\it Proof.}   

$i)$ One has
$$
 q^{2k} \left(  \frac {p^h}{q^k} \right) ' =  h p' p^{h-1}q^k -  kp^h q' q^{k-1} =  p^{h-1} q^{k-1} ( h p' q -  k q' p )= 0
$$
Then the function $\dst{\frac {p(x)^h}{q(x)^k}}$ is constant on $I$, hence the thesis.

$ii)$ Since analytic functions have isolated zeroes, there exists an interval $I^* \subset I$ such that $p(x) \neq 0 \neq q(x)$ for all $x \in I^*$. Working as in point $i)$ one proves that there exist $c_p, c_q \in \R$, $c_p c_q \neq 0$, such that $ c_q \, q(x) ^k + c_p\,   p(x)^h = 0$   for all $x \in I^*$. Then by analitycity such an inequality holds on all of $I$. 

$iii)$ In both $i)$ and $ii)$ the equality (\ref{equalemmauno}) becomes
$$
h p'(x) q(x) - hl q'(x) p(x) = h \Big( p'(x) q(x) - l q'(x) p(x) \Big) = 0,
$$
hence $p'(x) q(x) - l q'(x) p(x) = 0 $  that gives the thesis.

\hfill  $\clubsuit$  \\

\begin{remark} The Lemma \ref{lemmadue}, point $i)$, does not hold if  $p(x)$ and $q(x)$ vanish in $I$. In fact, let us choose $h=k=1$,  $I = (-1,1)$ and let $\phi \in C^\infty (I,\R)$ be a function flat at 0, vanishing at 0  and positive for $x \neq 0$.  Consider the functions so defined
$$
p(x) = \phi(x), \quad x \in I, \qquad q(x) =\left\{ \matrix{ \phi(x), \qquad x \geq 0 \cr - \phi(x), \quad x <0 }\right. . 
$$
One has $p,q \in C^\infty(I,\R)$ and the equation (\ref{equalemmauno}) holds on all of $I$, but they are not proportional on all of $I$. In fact, their ratio is $-1$ on $(-1,0)$ and $1$ on $(0,1)$.
\end{remark}

Most of our results are proved in a strip $\Sigma$, rather than in $\R^2$, also in order to have the possibility to split the domain of $F(x,y)$ into strips where it is injective, even if not globally injective. 

We say that a function $q \in C^1(I,\R)$ has a simple zero in $x_0$ if $q(x_0)=0 \neq q'(x_0)$. If $q(x_0)= q'(x_0)=0$ we say that $x_0$ is  a double zero of $q$.

\begin{lemma}  \label{lemmatre} The following statements hold.
\begin{itemize}
\item[i)] A non-singular map $F\in C^1(\Sigma,\R^2)$ of type $(0,1)$  is invertible. Moreover $q_1$ does not vanish in $I$.  
\item[ii)] Let $F(x,y) = (p_1(x) y + p_0(x),Q(x,y)) \in C^1(\Sigma,\R^2)$ be a non-singular map.  If  $p_1$ does not vanish in $I$, then $F(x,y) $  is injective. 
\item[iii)] Let $F(x,y) = (p_1(x) y + p_0(x),Q(x,y))  \in C^2(\Sigma,\R^2)$ be a non-singular map.  Then  $p_1 $ has no simple zeroes. 
\end{itemize}
\end{lemma}
{\it Proof.}

 { i)}  One has $F(x,y) = \left( p_0(x) ,  q_1(x) y + q_0(x) \right) $.
Its Jacobian determinant is
$$
d_F = p_0' q_1 \neq 0.
$$
hence both $p_0' $ and $q_1 $ do not vanish. As a consequence $p_0$ is strictly monotone and invertible, hence one has
$$
F^{-1}(u,v) = \left(   p_0^{-1}(u), \frac {v - q_0(p_0^{-1}(u) }{q_1(p_0^{-1}(u))}\right).
$$

 { ii)} If $p_1(x) \neq 0$ in $I$ then every level set of $P(x,y)$ is connected, hence $F(x,y)$ is injective.

 { iii)} 
The Jacobian matrix of $F$ is
\begin{equation}  \label{jgm-m} 
 J_{F}(x,y) = \left( \matrix{   
  p_1' y + p_0' &     p_1 \cr Q_x&   Q_y}   \right) .
\end{equation}

If  $p_1(x_0) = 0$ for some $x_0 \in I$, then the Jacobian determinant of $J_{F}(x_0,y)$ is
\begin{equation}  \label{lemdouble} 
 d_{J_{F}}(x_0,y) =  - Q_y (x_0,y) \, \Big( p_1'(x_0)  y + {p}_0'(x_0) \Big) .
 \end{equation}
  If $p_1'(x_0) \neq 0$, the  $y$-polynomial $ p_1'(x_0)  y + {p}_0'(x_0) $ has non-vanishing leading coefficient, hence it vanishes at some point $y_0$, contradicting 
$d_{J_{F}}(x_0,y_0) \neq 0$. Hence $p_1'(x_0) = 0 = p_1(x_0) $, so that $x_0$ is a double zero of $p_1(x)$.
\hfill  $\clubsuit$  \\

The hypotheses of Lemma \ref{lemmatre}, $ii)$, do not imply surjectivity, even if $\Sigma = \R^2$. The map
$$
F(x,y) = \left( e^x y + \sin x, e^x y + \sin x + e^{-x} \right) .
$$
is a Jacobian map of type $(1,1)$ defined on $\R^2$, injective by Lemma \ref{lemmauno}, but  not surjective. In fact, setting $(u,v) =  \left( e^x y + \sin x, e^x y + \sin x + e^{-x} \right)$, one has $v = u + e^{-x} > u$, hence the half-plane $v \leq u$  contains no images of $F(x,y)$. 

\begin{theorem}  \label{teononsingF1-1} 
A non-singular map $F\in C^1(\Sigma,\R^2)$ of type $(1,1)$  is the composition of a quasi-triangular map and a triangular map.  If $F$ is polynomial and $F(0,0)=0$, then it is the composition fo two triangular maps, hence it has a polynomial inverse.
\end{theorem}
{\it Proof.}  
 If one among $p_1, q_1$ is identically zero  then lemma \ref{lemmatre}, point i), applies, hence we assume both not to vanish identically on $I$. The Jacobian determinant of $F(x,y) = \left( p_1y + p_0, q_1y + q_0  \right)$ is
$$
d_F(x,y) = ( p_1' q_1 - p_1 q_1')y + p_0' q_1 - p_1 q_0'.
$$ 
If by absurd $p_1' q_1 - p_1 q_1'$ does not vanish identically there exists $x_0 \in I$ such that $d_F(x_0,y)$ is a polynomial of degree 1 in $y$. Hence $d_F(x_0,y)$ vanishes at some $y_0$, contradicting the non-singularity of $F(x,y)$. As a consequence $p_1$ and $q_1$ satisfy the differential relation
\begin{equation}  \label{R1-1} 
p_1' q_1 - q_1' p_1 = 0
\end{equation}
on all of $I$.
The functions $p_1$ and $q_1$ do not vanish simultaneously, since $d_F(x,y) \neq 0$. Let  $x_0 \in \R$ be such that $q_1(x_0) \neq 0$. 
Let $I_0$ be the largest interval containing $x_0$ such that $q_1(x) \neq 0 $ on $I_0$. 
By Lemma \ref{lemmadue} there exist $c \in \R$ such that $ p(x) = c\, q(x)$   for all $x \in I_0$.
Let us assume by absurd that $I_0 \neq I$. Let $x_\partial$ be a point of the boundary $\partial \left( I_0 \setminus  I \right)$. By definition of $I_0$ one has   $q_1(x_\partial) = 0$.  Since $ p(x) = c\, q(x)$ on $I_0$, by continuity  one has $p_1(x_\partial) = c\, q_1(x_\partial) = 0$, that implies $d_F(x_\partial,y) = 0$, contradicting the non-singularity of $F(x,y)$. 
This proves that both $p_1$ and $q_1$ do not vanish in $I$, $ p(x) = c\, q(x)$   for all $x \in I$. Let us consider the triangular map $\Lambda(u,v) = (u-cv,v)$. The composed map
$$
\Gamma = \Lambda \circ F = \Big( p_1 y + p_0 - c q_1 y - c q_0, q_1 y -  q_0 \Big) =
\Big( p_0 - c q_0, q_1 y -  q_0 \Big) = \Big( \overline{p}_0 , q_1 y -  q_0 \Big) ,
$$
has Jacobian determinant $d_F$, since the determinant of $\Lambda$ is 1. It is a non-singular quasi-triangular map, hence it is invertible.  One has
$$
F = \Lambda^{-1} \circ \Gamma,
$$
that proves the thesis, since $ \Lambda^{-1}$ is triangular.

If $F(x,y)$ is polynomial, then one has $d_F = \overline{p}_0' q_1 \neq 0$, which is possible only if both $ \overline{p}_0' $ and $q_1 $ are non-zero constants. If $F(0,0)=0$, one has $ \overline{p}_0 = a  x $ and $q_1 = b$  for some $a,b \neq 0$.
Hence $\Gamma = (ax, by + q_0)$ is triangular and $F = \Lambda^{-1} \circ \Gamma$ is the composition of two triangular maps.

\hfill  $\clubsuit$  \\

Non-singular polynomial maps of type (1,1) can be easily found, for instance $F(x,y) = \left( (x^2  +1) y + 2x, (x^2  +1) y + x \right) $, whose Jacobian determinant is $d_F(x,y) = x^2 + 1$.  On the other hand, for a polynomial map of type (1,1) to be a Jacobian map it is necessary that both $p_1$ and $q_1$ be constant. In fact, if $d_F(x,y) \equiv d_F \in \R$, $d_F \neq  0$, from Lemma \ref{lemmadue} one has $q_1 = k_1 p_1$ for some $k_1 \neq 0$ and
$$
d_F= p_0' q_1 - p_1 q_0' = p_1 ( k_1 p_0' - q_0') . 
$$
This  is possible only if both $p_1$ and $k_1 p_0' - q_0'$ are non-zero constants. In such a case also $q_1 = k_1 p_1$ is constant. 

In next theorem, we set 
$$
d_{1m}(x) =  \Big| \matrix{ p_1(x)  &  q_1(x)  \cr p_m(x) & q_m (x)}  \Big| = p_1(x)   q_m (x) - q_1 (x) p_m(x)    ,
$$
$$
d_{1m}^*(x) =  \Big| \matrix{ p_1'(x)  &  q_1'(x)  \cr p_m(x) & q_m (x)}  \Big| = p_1'(x)   q_m (x) -  q_1' (x)p_m(x)   .
$$

\begin{theorem}  \label{teononsingFm-m} 
Let a non-singular map $F\in C^\omega(\Sigma,\R^2)$ be of the form
\begin{equation}  \label{} 
F(x,y) = \Big(  p_m(x) y^m + p_1(x) y + p_0(x) , q_m(x) y^m + q_1(x) y + q_0(x)  \Big) ,
\end{equation}
with $m >  1$, $p_m,q_m \in C^\omega(I,\R)$.  If one of the following condition holds,
\begin{itemize}
\item[i)] $d_{1m}(x) \neq 0$ for all $x \in J$,
\item[ii)] $d_{1m}(x) = 0 \  \Rightarrow   \ d_{1m}^*(x) \neq 0$ ,
\item[iii)] $m$ is even and  $d_{1m}(x) = 0 \  \Rightarrow   q_m(x) \neq 0  $ \ \  \Big($p_m(x) \neq 0  $ \Big),
\item[iv)] $m$ is odd  and  $d_{1m}(x) = 0 \  \Rightarrow   q_m(x) \neq 0 , \ q_1(x)  \leq 0  $ \ \ \Big($p_m(x) \neq 0,\ p_1(x) \leq 0  $ \Big).
\end{itemize}
then  $F$ is injective.  
\end{theorem}
{\it Proof.}  
One has
\begin{equation}  \label{jfm-m} 
 d_F(x,y) = \left| \matrix{ p_m' y^m + p_1' y + p_0' &  mp_m y^{m-1} + p_1 \cr q_m' y^m + q_1' y + q_0' &  m q_m y^{m-1} + q_1 }  \right| =
\end{equation}
$$
= m(p_m' q_m - q_m' p_m) y^{2m-1} + \dots +  p_0' q_1 - q_0' p_1 .
$$
Since $d_F(x,y)$ is an odd-degree polynomial in $y$, for $d_F(x,y)$ not to vanish, $p_m' q_m - q_m' p_m$ has to vanish identically. By Lemma  \ref{lemmadue} , point ii), there exists $k_m \neq 0$ such that  one has $p_m = k_m\, q_m$ in $I$. 
Let us consider the map $\dst{ \Lambda_{m}(u,v) = (u - k_m v,v)  }$. The composed map $G_{m} = \Lambda_{m} \circ F $ has determinant $d_F(x,y) \neq 0$ and has the form:
$$
G_m = \Lambda_{m} \circ F =  \Big( p_{m} y^{m} + p_1 y + p_0 - k_m ( q_m y^m +  q_1 y + q_0), q_m y^m +  q_1 y + q_0 \Big) =
$$    $$
= \Big( (p_{m} - k_m q_m^{m}) y^{m} +  (p_{1} - k_m q_1) y + p_{0} - k_m q_0  ,  q_m y^m + q_1 y + q_0 \Big) =
$$  $$
= \Big( \overline{p}_1 y + \overline{p}_0 ,   q_m y^m + q_1 y + q_0 \Big) .
$$
If $\overline{p}_1(x) = p_1(x)   - k_m   q_1(x)  \neq 0$ for all $x\in I$, by Lemma \ref{lemmauno} the level sets of $ \overline{p}_1 y + \overline{p}_0 $ are connected, that implies the injectivity of $G_m(x,y)$, hence that of $F(x,y)$.  

  { i)} Since $d_{1m} =  p_1   q_m  - p_m   q_1  =  p_1   q_m  - k_m q_m   q_1  = q_m \big(p_1   - k_m   q_1    \big) = q_m \overline{p}_1  \neq 0 $, one has $\overline{p}_1(x)  \neq 0$ for all $x\in I$, proving the injectivity of $G_m$, hence that of $F$.

  { ii)}  Let us consider the Jacobian matrix of $G_m$,
\begin{equation}  \label{jgm-m} 
 J_{G_m}(x,y) = \left( \matrix{ \overline{p}_1' y + \overline{p}_0' &   \overline{p}_1  \cr
 q_m' y^m + q_1' y + q_0' &  m q_m^{m-1} y^{m-1} + q_1}   \right) .
\end{equation}

If there exists $x_0 \in \R$ such that 
$$
\overline{p}_1(x_0) =  p_{1}(x_0) - k_m q_1(x_0)  =  0 ,
$$ 
then the Jacobian determinant of $J_{G_m}(x_0,y)$ is
\begin{equation}  \label{detjgm-m} 
 d_{J_{G_m}}(x_0,y) =  \Big(\overline{p}_1'(x_0) y + \overline{p}_0'(x_0) \Big) \Big(  m q_m^{m-1}(x_0) y^{m-1} + q_1 (x_0) \Big)  .
\end{equation}

By absurd, if $\overline{p}_1(x_0) = 0$ then $d_{1m}(x_0) = 0$, hence $d_{1m}^*(x_0) \neq 0$. Since $p_m$ and $q_m$ vanish at the same points, this implies  $p_m(x_0) \neq 0 \neq q_m(x_0)$. For all $x \in J$ such that $q_m(x) \neq 0$ one has
$$
\overline{p}_1'(x) = p_1' (x)  - k_m   q_1'(x) =  p_1' (x)  - \frac{p_m(x)}{q_m(x)}   q_1'(x) =
\frac{d_{1m}^*(x)}{q_m(x)}   
$$
hence
$$
\overline{p}_1'(x_0) =\frac{d_{1m}^*(x_0)}{q_m(x_0)}   \neq 0  ,
$$
contradicting Lemma  \ref{lemmatre}, iii).  
This proves that $\overline{p}_1(x)  \neq 0$ for all $x\in I$, again proving the injectivity of $G_m$, hence that of $F$.

 { iii)} Also in this case one has a contradiction assuming $\overline{p}_1(x_0) = 0$ for some $x_0$, because  $d_{1m} (x_0)=  q_m(x_0) \overline{p}_1 (x_0) = 0 $ and by assumption $q_m(x_0) \neq 0$.
Since $m-1$ is odd one has 
\begin{equation}  
\label{radice} 
  d_{J_{G_m}} \left( x_0, \sqrt[m-1] {-\frac { q_1(x_0)}{m  q_m^{m-1}(x_0)} }\right) = \Big(\overline{p}_1'(x_0) y + \overline{p}_0'(x_0) \Big)  \cdot 0 =0,
\end{equation}  
contradicting the non-singularity of $G_m$.

 { iv)} In this case  one has a contradiction because   $q_m(x_0) \neq 0$ and the term under the root of formula (\ref{radice}) is non-negative. \\

\hfill  $\clubsuit$  \\

If $p_m(x)$ or $q_m(x)$ do not vanish on $I$, point iii) of Theorem 1 applies. In general the condition $p_m(x) \neq 0$ does not imply the connectedness of the level sets of a non-singular analytic $P(x,y)$, as shown in  remark \ref{nonconnected} for a $y$-quadratic map. A special case is that of non-singular $y$-quadratic maps with non-vanishing leading coefficient.

\begin{corollary}  \label{cornonsingF2-2} 
A $y$-quadratic non-singular map $F\in C^\omega(\Sigma,\R^2)$ 
\begin{equation}  \label{} 
F(x,y) = \Big(  p_2(x) y^2 + p_1(x) y + p_0(x) , q_2(x) y^2 + q_1(x) y + q_0(x)  \Big) .
\end{equation}such that either $p_2(x) \neq 0$ or $q_2(x) \neq 0$ on $I$  is injective.  
\end{corollary}
{\it Proof.}  Point iii) of theorem  \ref{teononsingFm-m}  applies.

\hfill  $\clubsuit$  \\

The non-vanishing of the leading  coefficients of $F(x,y)$ is not necessary for the map non-singularity or injectivity. In fact, the map $F(x,y) = (x^3 y^2 + x ,x^3 y^2 + x + y)$ has Jacobian determinant $J_F(x,y) = 1 + 3 x^2 y^2 > 0$ and is injective, since it is the composition of the two maps  $F_1(x,y) = (x^3 y^2 + x ,y)$  and $F_2(u,v) = (u , u + v)$, both injective.

In next theorem we deal with another class of non-singular maps.

\begin{theorem}  \label{teononsingFL2h} 
Let $F\in C^\omega(\Sigma,\R^2)$ be a non-singular map of the form
\begin{equation}
F(x,y) = \left( \sum_{l=1}^L p_{2hl}y^{2hl} + p_1 y + p_0 , q_{2h}y^{2h} + q_0 \right) ,
\end{equation}
with $h > 0$  integer, $h \geq 1$. If $p_1(x) \neq 0$ on $I$, then $F(x,y)$ is injective
\end{theorem}
{\it Proof.}  By induction on $L$. For $L=1$ the map has the form
$$
F(x,y) = \left( p_{2h}y^{2h} + p_1 y + p_0 , q_{2h}y^{2h} + q_0 \right) .
$$
Applying the reduction procedure exposed at the beginning of theorem \ref{teononsingFm-m}, one reduces the injectivity of $F(x,y)$ to that of the map
$$
F(x,y) = \left(  p_1 y + \overline{p}_0 , q_{2h}y^{2h} + q_0 \right) .
$$
The coefficient $p_1(x)$ is not affected by such an operation - in other words, $\overline{p}_1(x) = p_1(x)$  - because in $Q(x,y)$ there are no terms of $y$-degree 1. 
By hypothesis one has $p_1(x) \neq 0$ on $I$, hence by lemma  \ref{lemmatre}, point ii), the map is injective.

Assuming the statement true for $L-1$, let us prove it for $L$. One has 
$$
d_F = \left(  {2h }p_{2hL}' q_{2h} - {2hL} p_{{2hL}} q_{2h}' \right) y^{2hL + 2h - 1} + \dots .
$$
$d_F(x,y)$ is a polynomial of  odd $y$-degree, hence the coefficient $ {2h }p_{2hL}' q_{2h} - {2hL} p_{{2hL}} q_{2h}'$ vanishes identically, hence $ p_{2hL}' q_{2h} - {L} p_{{2hL}} q_{2h}'$ on $I$. By lemma \ref{lemmadue} there exists $c \neq 0$ such that $p_{2hL} = c q_{2h}^L$.  Let us consider the triangular map $\dst{ \Gamma_{L}(u,v) = (u - c v^L, v)  }$. The composed map $G_{L} = \Gamma_{L} \circ F $ has determinant $d_F(x,y) \neq 0$ and has the form:
$$
\Gamma_{L}  \circ F =  \Big( \sum_{l=1}^L p_{2hl}y^{2hl} + p_1 y + p_0 - c \big(  q_{2h}y^{2h} + q_0  \big)^L, q_{2h}y^{2h} + q_0 \Big) =
$$  $$
 = \Big( p_{2hL} - c q_{2h}^{L} +\sum_{l=1}^{L-1} \overline{p}_{2hl}y^{2hl} + p_1 y + \overline{p}_0 ,  q_{2h}y^{2h} + q_0 \Big) = 
 $$  $$
 = \Big(  \sum_{l=1}^{L-1} \overline{p}_{lm}y^{lm} + p_1 y + \overline{p}_0  , q_m y^m + q_0 \Big) .
$$
Also in this case the coefficient $p_1(x)$ is not affected by by such an operation, hence the map $\Gamma_{L}  \circ F $ satisfies the induction hypotheses. This completes the proof.

\hfill  $\clubsuit$  \\

\section{$y$-polynomial maps with Jacobian determinant independent of $y$}

We say that $F \in C^1 ( \Sigma,\R)$ is a \textit{$\delta$-map} if there exists a map $\delta: I \to \R$ such that $d_F(x,y) = \delta(x)$ for all $x \in I$. For instance, every map of the form $F(x,y) = ( x,\psi(x) y)$ is a $\delta$-map, non-singular if  $\psi(x) \neq 0$. Jacobian maps are $\delta$-maps, with $\delta(x) \equiv const. \neq 0$.

\begin{theorem}  \label{teodelFm-1} 
Let $F\in C^\omega(\Sigma,\R^2)$ be a non-singular $\delta$-map of type $(m,1)$. Then $F(x,y)$   is the  composition of a quasi-triangular map and $m$   triangular maps. As a consequence, its inverse is the  composition of a quasi-triangular map and $m$   triangular maps.
If $F$ is polynomial and $F(0,0) =(0,0)$, then it is the composition of $m+1$ triangular maps, hence it has a polynomial inverse. 
\end{theorem}
{\it Proof.}  The map has the form
\begin{equation}  \label{Fm-1} 
F(x,y) = \Big( p_m(x) y^m + \dots + p_1(x) y + p_0(x) ,  q_1(x) y + q_0(x) \Big)  .
\end{equation}
Let us assume $d_F(x,y) = \delta(x) \neq 0$. 

If $m = 0$, one has $F(x,y) = \left( p_0(x) ,  q_1(x) y + q_0(x) \right) $, which is a quasi-triangular map,  injective  by Lemma \ref{lemmatre}, $i)$. If $F$ is polynomial, the relationship
$$
d_F(x,y) = p_0' q_1 \neq 0
$$
  is possible only if both $p_0'$ and $q_1$ are non-zero constants, hence $F$ is a triangular map. 

For $m > 0$ we prove the thesis by induction on $m$.
If $m = 1$, one has $F(x,y) = \left( p_1(x) y + p_0(x) ,  q_1(x) y + q_0(x) \right) $. One has
$$
 d_F(x,y) = \left| \matrix{  p_1' y + p_0' &    p_1 \cr   q_1' y + q_0' &   q_1 }  \right| = (p_1' q_1 - p_1 q_1')y + p_0'q_1 - p_1 q_0' = \delta \neq 0 .
$$
Since $\delta$ does not depend on $y$, one has $p_1' q_1 - p_1 q_1' = 0 $ for all $x \in I$.  By Lemma \ref{lemmadue} there exists $c_1 \neq 0$ such that $p_1(x) = c_1 q_1(x)$. The map $\Lambda_{1}(x,y) = (x - c_1 y,y)$ has Jacobian determinant 1 and inverse $\Lambda_{1}^{-1}(u,v) = (u + c_1 v,v)$. The composed map $G_{1} = \Lambda_1 \circ F $ has the form:
$$
G_{1} = \Lambda_{1} \circ F =  \Big(  p_1 y + p_0 - c_1 (q_1 y + q_0),  q_1 y + q_0 \Big) = 
\Big(   p_0 - c_1  q_0,  q_1 y + q_0 \Big) .
$$    
$G_1$ is a quasi-triangular map. One has
\begin{equation}   
F(x,y) =   \Lambda_{1}^{-1} \circ G_1,
\end{equation}  
where $\Lambda_{1}^{-1} $ is a triangular map, so proving the statement for $m = 1$. If $F$ is polynomial, then 
$$
d_F = (p_0' - c_1  q_0') q_1 \neq 0,
$$
hence both $p_0' - c_1  q_0'$ and $q_1$ are non-zero constants, say $p_0' - c_1  q_0' = a \in \R$ and $q_1 = b \in \R$. As a consequence $p_0 - c_1  q_0 = ax$ (because $p_0(0) = q_0(0) = 0$) and $G_1 = (a x ,by + q_0)$ is a triangular map.

Let us assume the thesis to be true for $m$, with $m \geq 1$, so that for every Jacobian map $F(x,y)$ of type $(m,1)$, $m \geq 1$, there exist  a quasi-triangular map $\Lambda_0$ and  triangular maps $\Lambda_j$, $ j =1, \dots, m$, such that
\begin{equation}  \label{compotri}
F = \Lambda_{m} \circ \dots \circ \Lambda_{0} .
\end{equation}  
Let us prove that the same occurs for $m+1$. One has
\begin{equation} \label{d-F}
d_F = y^{m+1} \Big(p_{m+1}' q_1 - ({m+1}) p_{m+1} q_1' \Big) + \dots +  p_0' q_1 - q_0' p_1 = \delta \neq 0.
\end{equation}
Since $\delta$ does not depend on $y$, the coefficient of $y^{m+1}$ vanishes identically, hence by Lemma  \ref{lemmadue}  there exists $c_{m+1} \in \R$, $c_{m+1}\neq 0$, such that $p_{m+1}(x) = c_{m+1} \, q_1(x)^{m+1}$. Let us consider  the  map $\dst{ \Lambda_{m+1}(u,v) = (u - c_{m+1} v^{m+1},v)  }$, having as inverse  $\Lambda_{m+1}^{-1}(u,v) = (u + c_{m+1} v^{m+1},v)$. The composed map $G_{m+1} = \Lambda_{m+1} \circ F $ has the form:
$$
G_{m+1} = \Lambda_{m+1} \circ F =  \Big( p_{m+1} y^{m+1} + \dots + p_1 y + p_0 - c_{m+1} (q_1 y + q_0)^{m+1},  q_1 y + q_0 \Big) =
$$    $$
= \Big( (p_{m+1} - c_{m+1} q_1^{m+1}) y^{m+1} + \overline{p}_{m} y^{m}\dots + \overline{p}_1 y + \overline{p}_0 ,  q_1 y + q_0 \Big) =
$$  $$
= \Big( \overline{p}_{m} y^{m}\dots + \overline{p}_1 y + \overline{p}_0 ,  q_1 y + q_0 \Big) .
$$
One has $ \det J_{G_{m+1}} = d_F \neq 0$, since the Jacobian determinant of $\Lambda_{m+1}$ is 1. By induction hypothesis there exist triangular maps $\Lambda_j$, $ j =1, \dots, m$, and a quasi-triangular map $\Lambda_0$ such that
\begin{equation}
G_{m+1} = \Lambda_{m+1} \circ F = \Lambda_{m+1} \circ  \Lambda_{m} \dots \circ \Lambda_{0} .
\end{equation}
$\Lambda_{m+1}$ is a   triangular map, its inverse is itself a   triangular map, hence one has
$$
F = \Lambda_{m+1}^{-1} \circ G_{m+1} = \Lambda_{m+1}^{-1} \circ \Lambda_{m} \circ \dots \circ \Lambda_{0} .
$$
If $F(x,y)$ is polynomial and $F(0,0)=(0,0)$, then by induction hypothesis $\Lambda_0$ is triangular, hence $F(x,y)$ is the composition of $m+1$ triangular maps.

\hfill  $\clubsuit$  \\

\begin{remark}
The assumption $F(0,0) =(0,0)$ is not a restriction, since one can consider the map  $\overline{F}(x,y) = F(x,y) - F(x_0,y_0)$, which vanishes at $(x_0,y_0)$ and satisfies the same jacobian hypotheses as $F(x,y) $. The injectivity of $\overline{F}(x,y)$ is equivalent to that of $F(x,y)$.
\end{remark}

In the following corollary we consider a class of maps with both components having $y$-degree greater than 1.

\begin{corollary}  \label{cordelFm-m} 
Let $F\in C^\omega(\Sigma,\R^2)$ be a $\delta$-map of the form
\begin{equation}  \label{Fm-m} 
F(x,y) = \Big( p_m(x) y^m + p_1(x) y + p_0(x) , q_m(x) y^m + q_1(x) y + q_0(x)  \Big) .
\end{equation}  
Then $F$ is the  composition of a quasi-triangular map and $m+1$   triangular maps, hence it is injective. If $F$ is polynomial and $F(0,0) =(0,0)$, then it is the composition of $m+2$ triangular maps. 
\end{corollary}
{\it Proof.}  
If one among $p_m$ and $q_m$ vanishes identically, then Theorem \ref{teodelFm-1} applies. If both functions do not vanish identically, one has
\begin{equation}  \label{jfm-m} 
 d_F(x,y) = \left| \matrix{ p_m' y^m + p_1' y + p_0' &  mp_m y^{m-1} + p_1 \cr q_m' y^m + q_1' y + q_0' &  m q_m y^{m-1} + q_1 }  \right| =
\end{equation}
$$
= m(p_m' q_m - q_m' p_m) y^{2m-1} + \dots +  p_0' q_1 - q_0' p_1 .
$$
Since $F(x,y)$ is a $\delta$-map, the coefficient of $y^{2m-1}$ vanishes identically, hence by Lemma  \ref{lemmadue}  there exists $c_m \in \R$, $c_m\neq 0$, such that $p_{m}(x) = c_m \, q_m(x)$. Let us consider  the  map $\dst{ \Gamma(u,v) = (u - c_m v,v)  }$, having as inverse  $\Gamma^{-1}(u,v) = (u + c_m v,v)$. The composed map $ \Gamma \circ F $ has the form:
$$
\Gamma \circ F =  \Big( p_m y^m +  p_1 y + p_0  - c_m ( q_m y^m +  q_1 y + q_0)  , q_{m} y^{m} + q_1 y + q_0\Big) =
$$    $$
= \Big(  (p_{m} - c_m q_m) y^{m} +  (p_{1} - c_m q_1) y + p_{0} - c_m q_0,  q_{m} y^{m} + q_1 y + q_0 \Big) =
$$  $$
= \Big( (p_{1} - c_m q_1) y + p_{0} - c_m q_0,  q_{m} y^{m} + q_1 y + q_0   \Big) .
$$
The map $\Gamma \circ F$ has the same determinant as $F$, hance it is a $\delta$-map of type $(1,m)$. By Theorem  \ref{teodelFm-1} there exists a quasi-triangular map $\Lambda_0$ and $m$ triangular maps $\Lambda_j$, $j=1, \dots , m$ such that 
$$
\Gamma \circ F =  \Gamma \circ \Lambda_m \circ \dots \circ \Lambda_1 \circ \Lambda_0.
$$
Then one has 
$$
 F =   \Gamma^{-1} \circ \Lambda_m \circ \dots \circ \Lambda_1 \circ \Lambda_0.
$$
The polynomial case can be treated as in Theorem \ref{teodelFm-1}.

\hfill  $\clubsuit$  \\

As a special case we have the following corollary, concerned with the case $m = 2$.

\begin{corollary}  \label{cordelF2-2} 
Let $F\in C^\omega(\Sigma,\R^2)$ be a $y$-quadratic  $\delta$-map.  Then $F$ is the  composition of a quasi-triangular map and $3$   triangular maps. If $F$ is polynomial and $F(0,0) =(0,0)$, then it is the composition of $4$ triangular maps. 
\end{corollary}

In particular one has the following statement.

\begin{corollary}  \label{corjacF2-2} 
Let $F\in C^\omega(\Sigma,\R^2)$ be a $y$-quadratic Jacobian map.  Then $F$ is the  composition of a quasi-triangular map and $3$   triangular maps. If $F$ is polynomial and $F(0,0) =(0,0)$, then it is the composition of $4$ triangular maps. 
\end{corollary}

The procedure exposed in the proof of theorem ... can be applied to prove a similar result for another class of maps.

\begin{theorem}  \label{teodelFLm-m} 
Let $F\in C^\omega(\Sigma,\R^2)$ be a non-singular $\delta$-map of the form
\begin{equation}
F(x,y) = \left( \sum_{l=1}^L p_{lm}y^{lm} + p_1 y + p_0 , q_m y^m + q_0 \right) ,
\end{equation}
with $m$  integer, $m \geq 2$. Then $F(x,y)$   is the  composition of a quasi-triangular map and $L + m$   triangular maps. 
 If $F$ is polynomial and $F(0,0) =(0,0)$, then it is the composition of $L+m+1$ triangular maps, hence it has a polynomial inverse. 
\end{theorem}
{\it Proof.}  By induction on $L$. For $L=1$ the statement is a special case of corollary  \ref{cordelFm-m}, with $q_1(x) = 0$ on $I$.  \\
Assuming the statement true for $L-1$, let us prove it for $L$. One has 
$$
d_F = \left(  mp_{Lm}' q_m - Lm p_{Lm} q_m' \right) y^{Lm + m - 1} + \dots - p_1 q_0'.
$$
$F(x,y)$ is a $\delta$-map, hence the coefficient $mp_{Lm}' q_m - Lm p_{Lm} q_m'   $ vanishes identically. By lemma  \ref{lemmadue} there exist $c \neq 0$ such that $  p_{Lm}  = c q_m^{L}$ on $I$.  
Let us consider the triangular map $\dst{ \Gamma_{L}(u,v) = (u - c v^L, v)  }$. The composed map $G_{L} = \Gamma_{L} \circ F $ has determinant $d_F(x,y) \neq 0$ and has the form:
$$
\Gamma_{L}  \circ F =  \Big( \sum_{l=1}^L p_{lm}y^{lm} + p_1 y + p_0 - c \big(   q_m y^m + q_0  \big)^L, q_m y^m + q_0 \Big) =
$$  $$
 = \Big( p_{Lm} - c q_m^{L} + \sum_{l=1}^{L-1} \overline{p}_{lm}y^{lm} + p_1 y + \overline{p}_0  , q_m y^m + q_0 \Big) = 
 $$  $$
 = \Big(  \sum_{l=1}^{L-1} \overline{p}_{lm}y^{lm} + p_1 y + \overline{p}_0  , q_m y^m + q_0 \Big) .
$$
The map $\Gamma_{L}  \circ F $ is a $\delta$-map and by induction hypothesis is the composition of a quasi-triangular map and $L - 1 + m$ triangular maps. Since $\Gamma_{L}^{-1} $ is a triangular map, the equality $F =    \Gamma_{L}^{-1}  \circ   \Gamma_{L}  \circ F$ proves the thesis.

The polynomial case can be treated as in theorem \ref{teodelFm-1}.
\hfill  $\clubsuit$  \\

\end{document}